# A BASIS FOR SLICING BIRKHOFF POLYTOPES

TREVOR GLYNN

ABSTRACT. We present a change of basis that may allow more efficient calculation of the volumes of Birkhoff polytopes using a slicing method. We construct the basis from a special set of square matrices. We explain how to construct this basis easily for any Birkhoff polytope, and give examples of its use. We also discuss possible directions for future work.

## 1. INTRODUCTION

The $n$th *Birkhoff polytope*, which we denote $B_n$, is the convex hull of the set of nonnegative $n \times n$ matrices whose rows and columns sum to 1. Much has been written about these polytopes and their properties. (See [1], [4], and [5], for example.) Of these properties, the volume is of special interest. The exact volume of $B_n$ is known only for $n$ up to 10. Beck and Pixton gave the exact volume of $B_{10}$ in [1]; the calculation took over a year using dozens of computers. In [3], Canfield and McKay gave an asymptotic formula for the volume of $B_n$. De Loera et. al. [7] found an exact formula for the volume of $B_n$, but its calculation remains intractable for large values of $n$.

Results from Liu [6] offer a new algorithm to calculate the volume of a polytope by slicing. The method requires the polytope to satisfy special integrality conditions. To describe these conditions, it is necessary to introduce some basic definitions. The following discussion is a simplified summary of the material covered in [6]. For the sake of simplicity, we will use the standard basis instead of an arbitrary one.

We concern ourselves with the real vector space $\mathbb{R}^D$ and the lattice $\mathbb{Z}^D$. First, we define two functions: Let $\pi^{(k)} \colon \mathbb{R}^D \to \mathbb{R}^k$ be the function that maps the vector $(x_1, x_2, \ldots, x_D)$ to the vector $(x_1, x_2, \ldots, x_k)$. Let $\pi_k \colon \mathbb{R}^k \to \mathbb{R}^D$ be the function that maps the vector $\vec{y}$ to its inverse image under $\pi^{(k)}$. In effect, one can think of $\pi_k(\vec{y})$ as shifting the entire affine space $\mathbb{R}^{D-k}$ by the vector $\vec{y}$. We will call the affine space $\pi_k(\vec{y})$ a *slicing space*, for reasons that will become clear shortly.

Now we are ready to introduce the special integrality conditions. We say that an affine space $U \subset \mathbb{R}^D$ is *integral* if
$$\pi^{(\dim U)}(U \cap \mathbb{Z}^D) = \mathbb{Z}^{\dim U}.$$
Likewise, we say a polytope is *affinely integral* if its affine hull is integral. Finally, we say that a polytope is *$k$-integral* if all of its faces of dimension less than or equal to $k$ are affinely integral.

The condition of $k$-integrality is very strong and difficult to satisfy for high values of $k$. Fortunately, there is a weaker condition that will suffice for our needs. We say that an affine space $U \subset \mathbb{R}^D$ is *in general position* if
$$\pi^{(\dim U)}(U) = \mathbb{R}^{\dim U}.$$

---

Trevor Glynn is partially supported by NSF grant DMS-1265702.

We define what it means for a polytope is in $k$-*general position* in the predictable way. A polytope is in *affinely general position* if its affine hull is in general position. A polytope is in $k$-general position if is faces of dimension less than or equal to $k$ are in affinely general position.

It is a theorem of [6] that for all $d$-dimensional polytopes $P \subset \mathbb{R}^D$, the following holds: If $P$ is $(k-1)$-integral and in $k$-general position, then the normalized volume of $P$, which we will write $\mathrm{NVol}(P)$, is

$$(1.1) \qquad \mathrm{NVol}(P) = \sum_{\vec{y} \in \pi^{(k)}(P) \cap \mathbb{Z}^k} \mathrm{Vol}_{D-k}\bigl(\pi_k(\vec{y}) \cap P\bigr),$$

where $\mathrm{Vol}_{D-k}(Q)$ is the volume of $Q$ with respect to the lattice $\mathbb{Z}^{D-k}$.

Algorithmically, one can imagine the formula in this way:
  (i) Project $P$ into $\mathbb{R}^k$.
  (ii) Collect the integral points of the projection.
  (iii) Derive the slicing spaces from the integral points, and intersect them with the original polytope $P$. This step creates "slices" of $P$, which will have dimension $D-k$ or less.
  (iv) Sum the volumes of the slices.

**Example 1.** In the case where $k = 1$, the slicing spaces are hyperplanes. Figure 1 depicts an example of this case. This is Example 2.1 in [6].

Let $P$ be the polytope in Figure 1. We can apply formula (1.1) to $P$ as follows:

$$\sum_{\vec{y} \in \pi^{(1)}(P) \cap \mathbb{Z}} \mathrm{Vol}\bigl(\pi_1(\vec{y}) \cap P\bigr) = \sum_{\vec{y} \in \{0,1,2,3,4\}} \mathrm{Vol}\bigl(\pi_1(\vec{y}) \cap P\bigr)$$
$$= 0 + 1 + 4 + 3 + 0$$
$$= 8$$

Note that the 1-dimensional polytopes $\{(0,0,0)\}$ and $\{(4,0,0)\}$ are still considered slices, although they contribute nothing to the total volume.

Notice how none of the edges of the polytope in Figure 1 lie parallel to the slices. This is an important observation; it demonstrates a general rule that is worth emphasizing.

**Lemma 2.** *A polytope is in $1$-general position if and only if none of its edges are parallel with one (equiv. any) of the slicing spaces.*

*Proof.* This follows from the definition of the $k$-general position property. □

Formula (1.1) only works on polytopes that are $(k-1)$-integral and in $k$-general position. Nevertheless, it is possible to apply the algorithm to *any* rational polytope by performing a special change of basis and then scaling the coordinates.

We plan to apply formula (1.1) to Birkhoff polytopes with $k = 1$. It is our hope that this may enable the derivation of a more tractable formula for the volume of $B_n$. To do this, we need to ensure the Birkhoff polytopes are 0-integral and in 1-general position. Fortunately, the property of being 0-integral is equivalent to being integral. Therefore one of the prerequisites is already satisfied.

However, it turns out that none of the Birkhoff polytopes beyond $B_2$ are in 1-general position. In this paper we devise a change of basis that will put any Birkhoff polytope in 1-general position. Using this change of basis, it will be possible to apply formula (1.1) to any Birkhoff polytope, with $k = 1$.

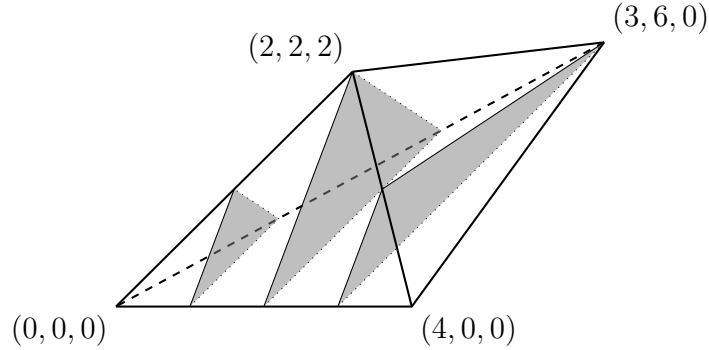

FIGURE 1. A polytope that is 1-integral, with its slices drawn and shaded.

Before we explain how we find this change of basis, it is worthwhile to explain why most of the Birkhoff polytopes are not in 1-general position. It is known that the direction vectors of the edges of the $n$th Birkhoff polytope comprise a special set of $(-1, 0, 1)$-matrices, which we will call $\mathcal{M}_n$. These matrices are discussed in the work of Brualdi and Gibson [2]. Here we give only a summary.

The matrices in $\mathcal{M}_n$ correspond with directed simple cycles in the complete bipartite graph $K_{n,n}$, with elements $u_1, u_2, \ldots, u_n$ and $v_1, v_2, \ldots, v_n$. The correspondence is as follows: A value of 1 for element $M[i, j]$ corresponds with an edge from $u_i$ to $v_j$, and a value of $-1$ corresponds with an edge *to $u_i$ from $v_j$*. In particular, if $n > 2$, then we can form a simple cycle in $K_{n,n}$ that includes neither $u_1$ nor $v_1$. Such a cycle will correspond with an edge in $B_n$ that is not in affinely general position. Hence $B_n$ is not in 1-general position for all $n > 2$.

## 2. MAIN RESULT

From Lemma 2, we know that we need to change to a basis such that none of the slicing spaces (in this case, hyperplanes) of $B_n$ are parallel to any edges. We will work backwards— first identifying a set of hyperplanes by their shared normal vector, and then deriving the basis from this set.

**Definition 3.** For all integers $n > 1$, let $\mathbf{V}_n = (a_{ij})$ be the square matrix

$$\begin{pmatrix} 0 & 0 & 0 & \cdots & 0 \\ 0 & 1 & 2 & \cdots & n-1 \\ 0 & n & 2n & \cdots & (n-1)n \\ \vdots & \vdots & \vdots & \ddots & \vdots \\ 0 & n^{n-2} & 2n^{n-2} & \cdots & (n-1)n^{n-2} \end{pmatrix},$$

and define $\mathbf{V}_1 = 1$.

We will prove the following theorem:

**Theorem 4.** *For all integers $n > 0$ and edges $e$ of $B_n$, the inner product $\mathbf{V}_n \cdot \mathrm{dir}\, e$ is nonzero.*

Thus we will show that the matrix $\mathbf{V}_n$ is not orthogonal to the direction vector $\mathrm{dir}\, e$ of any edge $e$ of the $n$th Birkhoff polytope. From this fact, we will conclude that any hyperplane normal to $\mathbf{V}_n$ satisfies the condition in Lemma 2. Before we prove this result, we will need to make several definitions.

First of all, since the matrices $\mathbf{M}$ represent simple cycles in $K_{n,n}$, there can never be more than one 1 or $-1$ in any row or column of $\mathbf{M}$. Furthermore, if a 1 (or $-1$) occurs in any particular row or column, then there must also be a $-1$ (resp. 1) in the same row or column. These observations will be useful later.

Our goal is to show that the dot product of $\mathbf{V}_n$ with any matrix $\mathbf{M} \in \mathcal{M}_n$ is never zero. We can visualize the product $\mathbf{V}_n \cdot \mathbf{M}$ as a special cycle consisting of elements in the matrix $\mathbf{V}_n$, where we consider two entries adjacent if they share a row or column. Since the matrix $\mathbf{M}$ contains negative entries, we must assign a sign to alternating elements in the cycle. This way, we do not need to bother with the set $\mathcal{M}$ at all, and can think just in terms of the matrix $\mathbf{V}_n$.

**Definition 5.** A *Birkhoff cycle* is a sequence of entries $a_{ij}$, $a_{ik}$, $a_{jk}$, ..., $a_{ij}$ of $\mathbf{V}_n$ that correspond with the entries in a matrix $\mathbf{M}$ in $\mathcal{M}_n$. We give these elements a sign according to the sign of the corresponding entry in $\mathbf{M}$. When we speak of the *sign* of an element in a Birkhoff cycle, it is this sign we are referring to, not the sign of the entry itself. (The vector $\mathbf{V}_n$ is nonnegative, so there should be no ambiguity.)

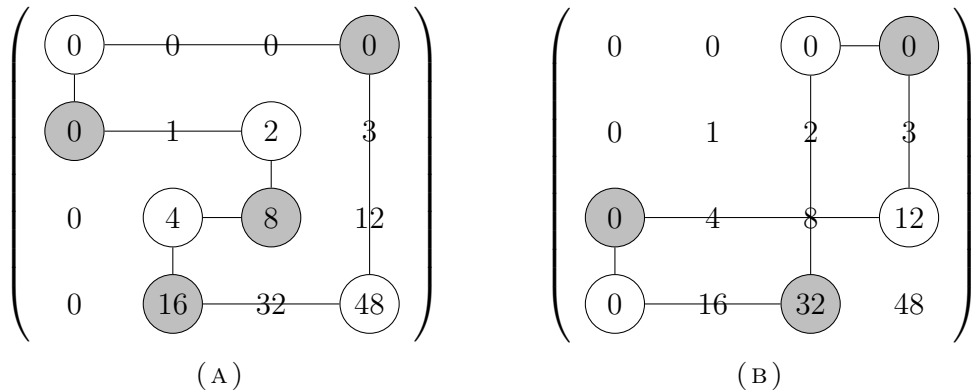

FIGURE 2. Two Birkhoff cycles in the matrix $\mathbf{V}_4$. The shaded elements are negative, and the unshaded ones are positive.

When we think of Birkhoff cycles, it can be useful to remember that their elements come from matrices. We may want to say that one element is "left of" or "above" another. When we say these things, we are referring to the elements' relative positions as matrix entries.

**Definition 6.** We say that an element $a_{rs}$ of a Birkhoff cycle in $\mathbf{V}_n$ is the *maximal element* if, for all other cycle elements $a_{ij}$, we have either $r > i$, or $r = i$ and $s > j$. Informally, the maximal element is the bottom and rightmost element.

**Example 7.** In the cycle in Figure 2a, the maximal element is 48. In the cycle in Figure 2b, the maximal element is 32.

Giving a sign to all of the elements in a Birkhoff cycle can be done in two ways, and it will be helpful for us to distinguish between the two.

**Definition 8.** If the maximal element is positive, we will say that the cycle itself is *positive*. Likewise, if the maximal element is negative, the cycle is *negative*.

**Example 9.** Figure 2a is a positive cycle, and Figure 2b is a negative cycle.

Since Birkhoff cycles are intended to represent the dot product $\mathbf{V}_n \cdot \mathbf{M}$, which is a scalar quantity, there should be a natural way of expressing that scalar as a property of the cycle itself. The following definition satisfies this requirement.

**Definition 10.** The *sum* of a Birkhoff cycle is the sum of its elements, signs included.

**Example 11.** The sum of the cycle in Figure 2a is 30, and that of the cycle in Figure 2b is $-20$.

*Notation.* If $c$ is a Birkhoff cycle, we will use $\sum c$ to denote its sum.

Notice that, in the previous example, the positive cycle has a positive sum, whereas the negative cycle has a negative sum. It turns out that this is not merely a coincidence, as we will prove now.

**Lemma 12.** *The sum of a positive Birkhoff cycle is positive. The sum of a negative Birkhoff cycle is negative.*

*Proof.* Since we can turn a negative cycle into a positive one by reversing all of the signs, we can restrict ourselves to the positive case. Let $c$ be an arbitrary positive cycle whose maximal element is in row $r$. Recall that $c$ must contain exactly zero or two nonzero elements from each row and column of $\mathbf{V}_n$. With this in mind, we know that the two cycle elements in row $r$ contribute a total value of no less than $n^{r-2}$ to $\sum c$. This happens then the maximal element is immediately to the right of its (negative) neighbor. Since every row of $\mathbf{V}_n$ is an arithmetic sequence, this is true no matter which column contains the maximal element.

Now we ask ourselves what the rows above $r$ contribute. We can overestimate by forgetting the positive elements and just taking the sum of the negative ones. Let $i_1 < i_2 < \ldots < i_s$ be the rows above $r$ with elements that are part of the cycle. We know from the definition of $\mathcal{M}$ that each of these rows contains exactly one element that is in $c$ and is negative. For each row $i_k$, let $j_k$ be the column containing this element. The columns $j_1, j_2, \ldots, j_s$ must be distinct because $c$ cannot contain two negative elements from any column of $\mathbf{V}_n$.

The sum of the negative elements above row $r$, not including their sign, is

(2.1) $$S = \sum_{k=1}^{s} (j_k - 1) n^{i_k - 2}.$$

In the right-hand side of equation (2.1), the terms that make up each summand are $j_k - 1$ and $n^{i_k - 2}$. Since the row indices $i_k$ are in increasing order, we know the terms $n^{i_k - 2}$ are also in increasing order. However, we know nothing about the order of the column indices $j_k$. The value of $S$ would be greatest if the column indices were also in increasing order. Therefore, if we let $j'_1 < j'_2 < \cdots < j'_s$ be the same column indices, then

$$\sum_{k=1}^{s} (j_k - 1) n^{i_k - 2} \leq \sum_{k=1}^{s} (j'_k - 1) n^{i_k - 2}.$$

We can go even further by choosing the maximum possible values of $i_k$ and $j'_k$, and maximizing their number.

$$\sum_{k=1}^{s}(j'_k - 1)n^{i_k-2} \leq \sum_{k=1}^{r-2}(n-r+k+1)n^{k-1} = \sum_{k=0}^{r-3}(n-r+k+2)n^k$$

Now we can say that

$$\sum c \geq n^{r-2} - S$$
$$\geq n^{r-2} - \sum_{k=0}^{r-3}(n-r+k+2)n^k.$$

When $r = 3$, we get the result $\sum c \geq 0$. In this case, we can make the inequality strict by simply checking that the "worst case" value of $S$ is unachievable unless there exists at least one positive cycle element above row $r$.

When $r > 3$, we get

$$\sum c \geq n^{r-2} - \sum_{k=0}^{r-3}(n-r+k+2)n^k$$
$$= n^{r-2} - \left((n-1)n^{r-3} + (n-2)n^{r-4} + \cdots + (n-r+2)\right)$$
$$= n^{r-2} - \left(n^{r-2} - n^{r-3} + n^{r-3} - 2n^{r-4} + \cdots + n - (r-2)\right)$$
$$= n^{r-4} + 2n^{r-5} + \cdots + (r-4)n + (r-2)$$
$$> 0.$$

We conclude that $c$ has a positive sum. Since we chose the cycle $c$ arbitrarily, it follows that all positive Birkhoff cycles have positive sums. □

Now we can give the proof of Theorem 4.

*Proof.* By applying Lemma 12, the proof of Theorem 4 is trivial: Every edge of $B_n$ corresponds with a matrix $\mathbf{M} \in \mathcal{M}$, and the dot product $\mathbf{V}_n \cdot \mathbf{M}$ is the sum of the Birkhoff cycle, which we have just shown to be nonzero. □

3. Constructing the Basis

To construct a new basis for the $n$th Birkhoff polytope, we begin by choosing $n^2 - 1$ linearly independent vectors from the hyperplane $H_n$ whose normal vector is parallel to $\mathbf{V}_n$. Since $\mathbf{V}_n$ is not orthogonal to any of the edge direction vectors of $B_n$, we know that $H_n$ is not *parallel* with any of the edges. Hence we can use this hyperplane to take slices of $B_n$.

There are several obvious choices for the particular vectors that will comprise our basis. For instance, there exist $2n - 1$ elementary $(0, 1)$-matrices with a single nonzero element in the top row or leftmost column. Since this row and column are filled with zeros in $\mathbf{V}_n$, clearly these elementary matrices lie in $H_n$.

For the rest of the basis, we take the $(n-1)^2 - 1$ vectors $v$ that are zero everywhere except for $v[2, 2] = -a_{ij}$ and $v[i, j] = 1$, where $i, j > 1$ and $i$ and $j$ are not both 2. The final vector will be the elementary matrix with a one in the second row and second column.

**Example 13.** The matrices in (3.1) form the basis for $\mathbf{V}_3$.

$$\text{(3.1)} \quad \begin{pmatrix} 1 & 0 & 0 \\ 0 & 0 & 0 \\ 0 & 0 & 0 \end{pmatrix} \begin{pmatrix} 0 & 1 & 0 \\ 0 & 0 & 0 \\ 0 & 0 & 0 \end{pmatrix} \begin{pmatrix} 0 & 0 & 1 \\ 0 & 0 & 0 \\ 0 & 0 & 0 \end{pmatrix} \begin{pmatrix} 0 & 0 & 0 \\ 1 & 0 & 0 \\ 0 & 0 & 0 \end{pmatrix} \begin{pmatrix} 0 & 0 & 0 \\ 0 & 0 & 0 \\ 1 & 0 & 0 \end{pmatrix}$$
$$\begin{pmatrix} 0 & 0 & 0 \\ 0 & -2 & 1 \\ 0 & 0 & 0 \end{pmatrix} \begin{pmatrix} 0 & 0 & 0 \\ 0 & -3 & 0 \\ 0 & 1 & 0 \end{pmatrix} \begin{pmatrix} 0 & 0 & 0 \\ 0 & -6 & 0 \\ 0 & 0 & 1 \end{pmatrix} \begin{pmatrix} 0 & 0 & 0 \\ 0 & 1 & 0 \\ 0 & 0 & 0 \end{pmatrix}$$

Earlier we mentioned the possible need for scaling the coordinates of the polytope after changing basis. But it turns out that this basis is unimodular. (This is easy to see by making a matrix where the $n$th row is formed by concatenating the rows of the $n$th basis vector.) Therefore, there is actually no need to scale the coordinates of the polytope; changing basis is all that is required.

Using this change of basis on a Birkhoff polytope will put it in 1-general position, which will allow the use of the slicing method described in [6].

There is an elegant interpretation of this result. The first eight basis vectors in (3.1) span a hyperplane that is parallel to every slice of $B_3$. The ninth basis vector is the offset that separates the slices.

After applying the change of basis to the Birkhoff polytopes $B_3$ and $B_4$, the resulting polytopes have the vertex representations shown in Figure 3.

$B_3$: $\{(1,0,0,1,0,0,1,0,0), (2,0,1,0,0,1,1,0,0), (3,0,0,1,1,0,0,1,0),$
$(5,1,0,0,0,1,0,1,0), (6,0,1,0,1,0,0,0,1), (7,1,0,0,0,0,0,0,1)\}$

$B_4$: $\{(6,0,0,0,1,0,1,0,0,1,0,0,1,0,0,0), (7,0,0,1,0,0,0,1,0,1,0,0,1,0,0,0),$
$(9,0,0,0,1,0,0,0,0,0,1,0,1,0,0,0), (11,0,1,0,0,0,0,1,0,0,1,0,1,0,0,0),$
$(13,0,0,1,0,0,0,0,0,0,0,1,1,0,0,0), (14,0,1,0,0,0,1,0,0,0,0,1,1,0,0,0),$
$(18,0,0,0,1,0,1,0,1,0,0,0,0,1,0,0), (19,0,0,1,0,0,0,1,1,0,0,0,0,1,0,0),$
$(24,0,0,0,1,1,0,0,0,0,1,0,0,1,0,0), (27,1,0,0,0,0,0,1,0,0,1,0,0,1,0,0),$
$(28,0,0,1,0,1,0,0,0,0,1,0,1,0,0,0), (30,1,0,0,0,0,1,0,0,0,0,1,0,1,1,0,0),$
$(33,0,0,0,1,0,0,0,1,0,0,0,0,0,1,0), (35,0,1,0,0,0,0,1,1,0,0,0,0,0,1,0),$
$(36,0,0,0,1,1,0,0,0,1,0,0,0,0,1,0), (39,1,0,0,0,0,0,1,0,1,0,0,0,0,1,0),$
$(44,0,1,0,0,1,0,0,0,0,1,0,0,1,0), (45,1,0,0,0,0,0,0,0,0,1,0,0,1,0),$
$(49,0,0,1,0,0,0,1,0,0,0,0,0,0,1), (50,0,1,0,0,0,1,0,1,0,0,0,0,0,0,1),$
$(52,0,0,1,0,1,0,0,0,1,0,0,0,0,0,1), (54,1,0,0,0,0,1,0,0,1,0,0,0,0,0,1),$
$(56,0,1,0,0,1,0,0,0,0,1,0,0,0,0,1), (57,1,0,0,0,0,0,0,0,0,1,0,0,0,0,1)\}$

FIGURE 3. The vertices of $B_3$ and $B_4$ after performing the change of basis to put them into 1-general position.

## 4. Future Work

Although our basis makes it possible to apply the slicing method to any Birkhoff polytope, it is unclear how best to do so. It is possible that the slicing method will yield an exact formula that is easier to compute, or that has interesting properties in itself.

It may be possible to invent new matrices, similar to the matrices $\mathbf{V}_n$, which are not orthogonal to any higher-dimensional facets of $B_n$. With such matrices, one could derive a basis that permits lower-dimensional slices.

## References

[1] Beck, M. and Pixton, D., 2003: The Ehrhart polynomial of the Birkhoff polytope. *Discrete Comput. Geom.*, **30**, 623–637.
[2] Richard A. Brualdi and Peter M. Gibson, 1976: Convex polyhedra of doubly stochastic matrices, IV. *Linear Algebra and Appl.* **15**, no. 2, 153–172.
[3] E. Rodney Canfield and Brendan D. McKay, 2007: The asymptotic volume of the Birkhoff polytope. arXiv:0705.2422.
[4] Chan, C. S. and Robbins D. P., 1999: On the volume of the polytope of doubly-stochastic matrices. *Experiment. Math.*, **8**, no. 3, 291–300.
[5] Diaconis, P. and Gangolli A., 1995: Rectangular arrays with fixed margins. *IMA Series on Volumes in Mathematics and its Applications,* # 72, Springer-Verlag, 15–41.
[6] Fu Liu, 2009: Higher integrality conditions, volumes, and Ehrhart polynomials. arXiv:0911.2051.
[7] Jesús A. De Loera, Fu Liu and Ruriko Yoshida, 2009: A generating function for all semi-magic squares and the volume of the Birkhoff polytope. *J. Algebraic Combin,* **30,** no. 1, 113–139, doi:10.1007/s10801-008-0155-y, arXiv:math/0701866.